\documentclass[a4paper,12pt]{article}
\usepackage{amsfonts,amssymb,epsfig}
\usepackage[latin1]{inputenc}
\newtheorem{thm}{Theorem}[section]
\newtheorem{cor}{Corollary}[section]
\newtheorem{lem}{Lemma}[section]
\newtheorem{pro}{Proposition}[section]
\newtheorem{dfn}{Definition}[section]
\newtheorem{rmq}{Remark}[section]
\newtheorem{expl}{Example}[section]
\newtheorem{nota}{Notation}[section]

\def\dessous#1\sous#2{\mathrel{\mathop{\kern0pt#2}\limits_{#1}}}

\newcommand{\al}{\alpha}
\newcommand{\be}{\beta}
\newcommand{\p}{\partial}
\newcommand{\n}{\nabla}
\newcommand{\pr}{\prime}
\newcommand{\mbb}{\mathbb}
\newcommand{\mf}{\mathfrak}
\newcommand{\mc}{\mathcal}
\newcommand{\nb}{\nonumber}
\newcommand{\beq}{\begin{eqnarray}}
\newcommand{\eeq}{\end{eqnarray}}
\newcommand{\bpro}{\begin{pro}}
\newcommand{\epro}{\end{pro}}
\newcommand{\blem}{\begin{lem}}
\newcommand{\elem}{\end{lem}}
\newcommand{\bdfn}{\begin{dfn}}
\newcommand{\edfn}{\end{dfn}}
\newcommand{\bcor}{\begin{cor}}
\newcommand{\ecor}{\end{cor}}
\newcommand{\bthm}{\begin{thm}}
\newcommand{\ethm}{\end{thm}}
\newcommand{\bex}{\begin{expl}}
\newcommand{\eex}{\end{expl}}
\newcommand{\brmq}{\begin{rmq}}
\newcommand{\ermq}{\end{rmq}}
\newcommand{\bnota}{\begin{nota}}
\newcommand{\enota}{\end{nota}}
\begin{document}
\title{$g$-Natural metrics of constant sectional curvature on
tangent bundles}
\author{S. Degla\footnote{Institut de Math\'ematiques et de
Sciences Physiques, 01 BP 613 IMSP Porto-Novo, BENIN.}$^{\, }$
\footnote{sdegla@imsp-uac.org (corresponding
author)}\quad , \quad 
J.-P. Ezin $^{\ast\, }$\footnote{jp.ezin@imsp-uac.org}   \quad and \quad 
L. Todjihounde $^{\ast\, }$\footnote{leonardt@imsp-uac.org}}  
\date{ }
\maketitle

\abstract{Let $(M,g)$ be a Riemannian manifold and $G$ a $g$-natural metric
on its tangent bundle $TM$. In this paper we prove first that the space
$(TM,G)$ has constant sectional curvature if and only if it is flat, and then 
we give a characterization of flat $g$-natural metrics on tangent bundles.}
\vspace{0.25cm}

\vspace{0.25cm}

\section*{Introduction}
In \cite{AS1}, K.M.T. Abbassi and M. Sarih introduced the notion
of $g$-natural me\-trics on the tangent bundle $TM$ of a
Riemannian manifold $(M,g)\,$. A metric $G$ on $TM$ is called a
$g$-natural metric if it comes from $g$ by a first order natural
operator 
$
S_+^2T^* \leadsto  (S^2T^*)T\, ,
$
where $S_+^2T^*$ and $ (S^2T^*)T$ denote respectively the natural
bundle of Riemannian metrics and the natural bundle of 
$(0,2)$-tensor fields on the tangent bundles (cf. \cite{KMS} for
the definitions of natural bundles and operators and associated
notions). They gave a characterization of $g$-natural me\-trics
on $TM$ in terms of functions defined on $\mbb{R}^+$, and 
obtained a necessary and sufficient conditions for $g$-natural
metrics to be either nondegenerate or Riemannian. But they did 
not give an explicit expression for the inverse of nondegenerate
$g$-natural metrics although it is important to compute
some \\ geometrical analysis tools like the Ricci tensor,
the scalar curvature,the Laplace operator, etc ... .

Some geometrical properties could be inherited on the 
$g$-natural metrics from the basic metric $g$ and conversely.
In \cite{AS2} the authors proved that if a  tangent
bundle equipped with a $g$-natural metric $(TM,\, G)$ is of
constant sectional curvature then the same holds for $(M,g)$.
Furthermore, making some restrictions on the Riemannian 
$g$-natural metrics on $TM$, the same authors gave the
characterization of flat Riemannian $g$-natural metrics on $TM$\\
(cf. \cite{AS3}).

In this paper we prove that if $(M,g)$ is non flat, its tangent bundle $TM$ equipped with a $g$-natural metric $G$  has non
constant sectional curvature, and also that only flat $g$-natural
metrics are of constant sectional curvature. 
  In  the next  section \ref{S2} we give some  preliminaries and 
 some known results on $g$-natural metrics. In  the section 
 \ref{S3} we compute explicitly the inverse of any nondegenerate
 $g$-natural metric. In section \ref{S4} using this inverse 
 expression and Koszul's formula, we determine the Levi-Civita
 connection of any nondege\-nerate $g$-natural metric. Finally in
 section \ref{S5}, we show that the flat Riemannian
  $g$-natural metrics are the only $g$-natural metrics that have
  a constant sectional curvature, then we give a characterization
  of these metrics. 
\section{Preliminaries}\label{S2}
Let $(M,g)$ be a Riemannian manifold and $\n$ the Levi-Civita 
connection of $g$. Then the tangent space of $TM$ at any point
$(x,u)\in TM$ splits into the horizontal and vertical subspaces
with respect to $\n$ :
$$ 
T_{(x,u)}TM = H_{(x,u)}M\oplus V_{(x,u)}M\; .
$$

If $(x,u)\in TM$ is given then, for any vector $X\in T_xM$, there
exists a unique vector $X^h\in H_{(x,u)}M$ such that 
$\pi_\ast X^h = X$, where $\pi : TM\rightarrow M$ is the natural
projection.  $X^h$ denotes the \textit{horizontal lift} of $X$ 
at the point $(x,u)\in TM$. The \textit{vertical lift} of a
vector $X\in T_xM$  at $(x,u)\in TM$ is a vector
$X^v\in V_{(x,u)}M$ such that $X^v.(df)=X.f$, for all functions 
$f$ on $M$. Here we consider $1$-forms $df$ on $M$ as functions
on $TM$ (i.e. $(df)(x,u) = u.f$). Note that the map $X\rightarrow
X^h$ is an isomorphism between the vector spaces $T_xM$ and 
$H_{(x,u)}M$. Similarly, the map $X\rightarrow X^v$ is an
isomorphism between the vector spaces $T_xM$ and $V_{(x,u)}M$.
Obviously, each tangent vector $\tilde{Z}\in T_{(x,u)}TM$ can be
written in the form $\tilde{Z} = X^h+Y^v$, where $X,Y\in T_xM$
are uniquely determined vectors.

If $\varphi$ is a smooth function on $M$, then 
\beq
X^h(\varphi\circ \pi) = (X\varphi)\circ \pi \mbox{ and }
X^v(\varphi\circ \pi) = 0
\eeq
hold for every vector field $X$ on $M$.

A system of local coordinates $\left(U\, ;\, x_i,\; i=1,\cdots,m\right)$ in
$M$ induces on $TM$ a system of local coordinates 
$\left(\pi^{-1}(U)\, ;\, x_i,u^i,\; i=1,\cdots,m\right)$.\\ Let 
$X=\sum_{i=1}^m X^i \frac{\p}{\p x_i}$ be the local expression in $U$ of a 
vector field $X$ on $M$. Then, the horizontal lift $X^h$ and the vertical lift
$X^v$ of $X$ are given, with respect to the induced coordinates, by :
\beq
X^h & = &\sum_iX^i\frac{\partial}{\partial x_i}-
\sum_{i,j,k}\Gamma_{jk}^iu^jX^k\frac{\partial}{\partial u^i}\quad\mbox{ and }\\
X^v & = & \sum_{i}X^i\frac{\partial}{\partial u^i}\, ,
\eeq
where the $(\Gamma_{jk}^i)$ are the   Christoffel's symbols of $g$.

Next, we  introduce some notations which will be used to describe vectors
obtained from lifted vectors by basic operations on $TM$. Let $T$ be a tensor
field of type $(1,s)$ on $M$. If $X_1,X_2,\cdots,X_{s-1}\in T_xM,$ then 
$h\{T(X_1,\cdots, u,\cdots,X_{s-1})\}$ (respectively $v\{T(X_1,\cdots, u,
\cdots,X_{s-1})\}$) is a horizontal (respectively vertical) vector at
$(x,u)$ which is defined by the formula
$$
h\{T(X_1,\cdots, u,\cdots,X_{s-1})\}=\sum u^\lambda
\left(T(X_1,\cdots,\left(\frac{\p}{\p x_\lambda}\right)_x\,,\cdots,X_{s-1})\right)^h
$$
$$
(\mbox{ resp.}\quad 
v\{T(X_1,\cdots, u,\cdots,X_{s-1})\}=\sum u^\lambda
\left(T(X_1,\cdots,\left(\frac{\p}{\p x_\lambda}\right)_x\,,\cdots,X_{s-1})\right)^v
\, ) .
$$
 In particular, if $T$ is the identity tensor of type $(1,1)$, then we
 obtain the geodesic flow vector field at $(x,u)$,
 $\xi_{(x,u)} =\sum_{\lambda}u^{\lambda}\left(\frac{\p}{\p x_\lambda}\right)_{(x,u)}^h$, and the cano\-nical vertical vector at $(x,u)$,
 $\mc{U}_{(x,u)}=
 \sum_{\lambda}u^\lambda\left(\frac{\p}{\p x_\lambda}\right)_{(x,u)}^v$.
 
Moreover  $h\{T(X_1,\cdots, u,\cdots,u,\cdots,X_{s-t})\}$
and $v\{T(X_1,\cdots, u,\cdots,u,\cdots,X_{s-t})\}$ are defined by 
similar way. 

Also let us make the notations
\beq
h\{T(X_1,\cdots,X_s)\} =: T(X_1,\cdots,X_s)^h
\eeq
 and 
\beq
 v\{T(X_1,\cdots, X_s)\}=: T(X_1,\cdots,X_s)^v\; .
\eeq
Thus $h\{X\}=X^h$ and $v\{X\}=X^v$, for each vector field $X$ on $M$.

From the preceding quantities, one can define vector fields on $TU$ in the
following way: If $u=\sum_iu^i\left(\frac{\p}{\p x_i}\right)_x$ is a given point in $TU$ and $X_1,\cdots,X_{s-1}$ are vector fields
on $U$, then we denote by 
$$
h\{T(X_1,\cdots, u,\cdots,X_{s-1})\}\quad
(\mbox{respectively}\quad
v\{T(X_1,\cdots, u,\cdots,X_{s-1})\}) $$
the horizontal (respectively vertical) vector field on $TU$ 
defined by
$$
h\{T(X_1,\cdots, u,\cdots,X_{s-1})\}=\sum_\lambda u^\lambda
T(X_1,\cdots,\frac{\p}{\p x_\lambda}\,,\cdots,X_{s-1})^h
$$
$$
(\mbox{ resp.}\quad
v\{T(X_1,\cdots, u,\cdots,X_{s-1})\}=\sum_\lambda u^\lambda
T(X_1,\cdots,\frac{\p}{\p x_\lambda}\,,\cdots,X_{s-1})^v
\, ) .
$$
Moreover, for vector fields $X_1,\cdots,X_{s-t}$ on $U$, where 
$s\, ,\, t\in \mbb{N}^\ast\, (s>t) $,  the\\
vector fields $h\{T(X_1,\cdots, u,\cdots, u, \cdots,X_{s-t})\}$ and \\
$v\{T(X_1,\cdots, u,\cdots,u, \cdots,X_{s-t})\}$, on $TU$, are defined by
similar way.

The Riemannian curvature of $g$ is defined by
\beq
R(X,Y) = \left[\n_X\, ,\, \n_Y\right] -\n_{[X,Y]}\; .
\eeq
Now, for $(r,s)\in \mbb{N}^2$, we denote by 
$\pi_M:\ TM\rightarrow M$ the natural  projection and $F$ the natural bundle 
defined by  
\beq
FM&=&\pi_M^*(\underbrace{T^*\otimes \cdots\otimes T^*}_{\mbox{$r$ times}}\otimes
\underbrace{T\otimes \cdots\otimes T}_{\mbox{$s$ times}})M\rightarrow M, \\ \nb
Ff(X_x,S_x)&=&(Tf.X_x,(T^*\otimes \cdots\otimes T^*\otimes 
T\otimes \cdots\otimes T)f.S_x) 
\eeq
for all manifolds $M$, local diffeomorphisms $f$ of $M$, $X_x\in T_xM$ and\\
$S_x\in (T^*\otimes \cdots\otimes T^*\otimes
T\otimes \cdots\otimes T)_xM$. We call the sections of the canonical projection
$FM\rightarrow M$ $F$-tensor fields of type $(r,s)$.
  So, if $\oplus$ denotes the  fibered product of 
fibered  manifolds, then the  $F$-tensor fields are mappings\\  
$A:\ TM\oplus \underbrace{TM\oplus\cdots \oplus TM}_{\mbox{$s$ times}}
\rightarrow \sqcup _{x\in M}\otimes^rT_xM$ which are   linear in the last\\
$s$ summands and such that  $\pi_2\circ A=\pi_1$, where $\pi_1$ and
$\pi_2$ are respectively  the  natural projections of the source and target
fiber bundles of $A$.
For $r=0$ and $s=2$, we obtain the classical notion of $F$-metrics. So,
$F$-metrics are mappings $TM\oplus TM\oplus TM\rightarrow \mbb{R}$ 
which are linear in the second and  the third argument.

\bpro\cite{AS1}\label{P1}
Let $(M,g)$ be  a Riemannian manifold and $G$  a\\  $g$-natural metric on $TM$. 
Then if $dim\,M\geq 2$,   there exists six functions  $\al_i,\
 \be_i:\mbb{R}^+\to \mbb{R},\ i=1,2,3,$
such that for any $x\in M$ and all  vectors $u,\ X,\ Y\in T_xM$, we have
\beq
\left\{
\begin{array}{lcl}
G_{(x,u)}\left(X^h,Y^h\right)
   &=& (\al_1+\al_3)(t)g_x(X,Y)
   +(\be_1+\be_3)(t)g_x(X,u)g_x(Y,u),\\ \nb
   & & \\ \nb
G_{(x,u)}\left(X^h,Y^v\right)
   &=& \al_2(t)g_x(X,Y)
      +\be_2(t)g_x(X,u)g_x(Y,u),\\ \nb
       & & \\ \nb
G_{(x,u)}\left(X^v,Y^h\right)
   &=& \al_2(t)g_x(X,Y)
         +\be_2(t)g_x(X,u)g_x(Y,u),\\ \nb
  & & \\ \nb
G_{(x,u)}\left(X^v,Y^v\right)
   &=& \al_1(t)g_x(X,Y)
         +\be_1(t)g_x(X,u)g_x(Y,u), \nb
\end{array}
\right.
\eeq
where $t=g_x(u,u)$, $X^h$ and $X^v$ are  respectively the  horizontal lift
and the vertical lift  of the  vector $X\in T_xM$ at the  point  $(x,u)\in TM$.

For  $dim\,M=1$, the same holds with $\be_i=0,\quad i=1,2,3.$
\epro

\bnota\label{N1}\
\begin{itemize}
\item $\phi_i(t)=\al_i(t)+t\be_i(t)$,
\item $\al(t)=\al_1(t)(\al_1+\al_3)(t)-\al_2^2(t)$,
\item $\phi(t)=\phi_1(t)(\phi_1+\phi_3)(t)-\phi_2^2(t)$,
\end{itemize}
for all $t\in \mbb{R}^+$.
\enota

\bpro\cite{AS1}\label{P2}
 A $g$-natural metric $G$ on the tangent bundle of a\\ Riemannian  manifold
 $(M,\ g)$ is :
\begin{enumerate}
\item[(i)]  nondegenerate if and only if the functions $\al_i,\, \be_i,\,
i=1,2,3$ of  Proposition \ref{P1} defining $G$, satisfy
\beq
\al(t)\phi(t)\neq 0
\eeq
for all $t\in \mbb{R}^+$.
\item[(ii)]  Riemannian if and only if the functions  
$\al_i,\, \be_i,\,i=1,2,3$ of \\ Proposition \ref{P1}
defining $G$,  satisfy the  inequalities
\beq
\left\{
\begin{array}{ll}
\al_1(t)>0,& \phi_1(t)>0,\\
\al(t)> 0, &\phi(t)>0,
\end{array}
\right.
\eeq
for all $t\in \mbb{R}^+$.

For $dim\,M=1$, this system  reduces to $\al_1(t)>0$
and $\al(t)>0$, for all $t\in \mbb{R}^+$.
\end{enumerate}
\epro

The following lemmas will be useful in the sequel.
\blem\cite{D}\label{L1}
Let $(M,g)$ be a Riemannian manifold, $\n$ be the Levi-Civita
connection and $R$ be the Riemannian curvature of $g$. Then the
Lie bracket on the tangent bundle $TM$ of $M$ satisfies
\begin{enumerate}
\item $\left[X^h,Y^h\right]= \left[X,Y\right]^h
-v\left\{R(X,Y)u\right\}\, ,$
\item  $\left[X^h,Y^v\right]=\left(\n_XY\right)^v\, ,$
\item $\left[X^v,Y^v\right]=0\, ,$
\end{enumerate}
for all $X\,,\, Y\, , \, Z \in \mf{X}(M)$.
\elem
\blem\cite{AS2}\label{L1p}
Let $(M,g)$ be a Riemannian manifold, $(x,u)\in TM$ and 
$X,Y,Z\in \mathfrak{X}(M)$, $f$ a function defined from $\mbb{R}$
to $\mbb{R}$, and denote by $F_Y$  the function on $TM$ defined
by $F_Y(u)=g_x(Y_x,u)$, for all $(x,u)\in TM$. 
Then we have:
\begin{enumerate}
\item $X_{(x,u)}^h .f(|u|^2)=0$,
\item $X_{(x,u)}^v .f(|u|^2)=2f'(|u|^2)g_x(X_x,u)$,
\item $X_{(x,u)}^h . F_Y=g_x((\n_XY)_x,u)=F_{\n_XY}(x,u)$,
\item $X_{(x,u)}^h .(g(Y,Z)\circ \pi) = X_x.(g(Y,Z))$,
\item $X_{(x,u)}^v .(g(Y,Z)\circ \pi) = 0$,
\item $X_{(x,u)}^v .F_Y=g_x(X,Y)$,\\
where $|u|^2=g_x(u,u)\, $.

\end{enumerate}
\elem

From now on, whenever we consider an arbitrary Riemannian 
$g$-natural metric $G$ on $TM$, we implicitly assume that it is
defined  by the functions 
$\alpha_i, \beta_i:\mbb{R}^+\longrightarrow \mbb{R},\ i=1,2,3 $ 
given in   Proposition \ref{P1} .

All real functions $\alpha_i, \beta_i, 
\phi_i, \alpha, \mbox{ and }\phi$ and their derivatives are evaluated 
at $t:=g_x(u,u)\,$, $u\in T_xM$, unless otherwise stated.
    
\section{Inverse of nondegenerate  g-natural\\ metrics}\label{S3}
Let $(a,b)\in \mbb{R}^2,\ m\in \mbb{N}^\ast$,  
$u=(u^1,\ \cdots\ ,u^m)\in \mbb{R}^m$ 
and denote by  $\mu(a,b,u)$ the following square matrix of order  
$m\in \mbb{N}^\ast$ :
\beq
\mu(a,b, u)=
\left(
\begin{array}{ccc}
a+b(u^1)^2 &{}  & {}  \\
     {}         &   {}  &       bu^iu^j\\
               {}   &\ddots &       {}\\
	                 bu^iu^j &  {}  &{}\\  
{} &     {}      & a+b(u^m)^2
          \end{array}
 \right),
\eeq
that is $[\mu(a,b,u)]_{ij}=a\delta_{ij}+bu^iu^j\, $.

We establish the following lemma which is  easy to check by
straightforward computation:
\blem\label{L2}
 If $\quad a(a+b|u|^2)\neq 0$,  
then  $\mu(a,b,u)$ is invertible and its inverse 
$\mu(a,b,u)^{-1}$ is given by 
\beq
\mu(a,b,u)_{ij}^{-1} &=& 
\frac{\delta_{ij}}{a}-\frac{b}{a(a+b|u|^2)}u^iu^j,
\eeq
where $\mu(a,b,u)_{ij}^{-1}$ is the element of $i^{th}$ line and of $j^{th}$ 
column of the matrix $\mu(a,b,u)^{-1}\, $ and $|u|^2 =\sum_{i=1}^m (u^i)^2$.
\elem
 
Next, we are going to determine  the inverse of a nondegenerate $g$-natural
 metric $G$ .

Let $(U,  x_i ,i=1,\cdots,m)$ be a normal coordinates system of 
$(M,g)$ centred at $p\in M$, and 
$(\pi^{-1}(U); x_i,u^i,i=1,\cdots,m)$ its induced coordinates 
system on $TM$. For $l=1,2,3$; let us consider the matrix-value
functions
\newpage 
\beq
M_l(x,u)=\left(\,\al_lg_{ij}+\be_lg(\p_{x_i}\, ,u)g(\p_{x_j}\, ,u)\,
 \right)_{1\leq i,j\leq m }\; ,\forall\,  (x,u)\in \pi^{-1}(U)\, , 
\eeq
where  $\p_{x_i}=\frac{\p}{\p_{x_i}}$ and 
$g_{ij} = g(\p_{x_i},\p_{x_j})$ on $U$.\\
\newline
So $
\begin{array}{c}
\left(
\begin{array}{ccc}
(M_1 +M_3) & {}& M_2 \\
         &{} &    \\
M_2      & {}& M_1 \\      
\end{array}
\right)\, \\ 
{}
\end{array}
$ is the matrix-value functions of $G$ in the\\
local frame  $(\p_{x_i}^h,\p_{x_i}^v)_{i=1,\cdots,m}$ on
$\pi^{-1}(U)$  and we have
\beq
G\equiv \left(
\begin{array}{ccc}
(M_1 +M_3) & {}& M_2 \\
         &{} &    \\
M_2      & {}& M_1
\end{array}
\right)\, .
\eeq
If $G$ is nondegenerate, its inverse $G^{-1}$ has  the form 
\beq\label{e}
G^{-1}\equiv
\left(
\begin{array}{ccc}
\Lambda & {} & \Theta \\
  & {} &{} \\
  \Theta & {} & \Omega
  \end{array}
  \right)
\eeq
 where $ \Lambda =(\lambda^{ij})_{1\leq i,j\leq m}\, ,\ 
\Theta=(\theta^{ij})_{1\leq i,j\leq m}\, , \mbox{ and } \Omega\,
=(\omega^{ij})_{1\leq i,j\leq m} $  are square matrix-value
functions of order  $\,  m\, $, defined on $\pi^{-1}(U)$.

Therefore we have the following proposition:
\vspace{0.15 cm}

\bpro\label{P4} 
If 
\beq\label{Lp4}
\left\{
\begin{array}{l}
\al(t)\phi(t)\neq 0\\
\al_1(t)(\al_1+\al_3)(t)\neq 0\\
\phi_1(t)(\phi_1+\phi_3)(t)\neq 0
\end{array}
\right.
\eeq
for any $t\in \mbb{R}^+$, then the  blocks of the  matrix-value functions in (\ref{e}) satisfy :
\beq\label{1e}
\Lambda (p,u)&\equiv& \left(\lambda^{ij}(p,u)\right)_{1\leq i\, ,\, j\leq m} 
\mbox{ with } \\
\lambda^{ij}(p,u) &=& \frac{\alpha_1}{\alpha} \delta_{ij} 
-\psi_\lambda u^iu^j
\, ,
\eeq
\beq\label{2e}
\Theta(p,u) &\equiv& \left(\theta^{ij}(p,u)\right)_{1\leq i\, ,\, j\leq m}
\mbox{ with } \\
\theta^{ij}(p,u) &=&-\frac{\alpha_2}{\alpha} \delta_{ij} 
-\psi_\theta u^iu^j \, ,
\eeq
\beq\label{3e}
\Omega(p,u) &\equiv& \left(\omega^{ij}(p,u)\right)_{1\leq i\, ,\, j\leq m}
\mbox{ with } \\
\omega^{ij}(p,u) &=&\frac{\alpha_1+\alpha_3}{\alpha} \delta_{ij} 
-\psi_\omega  u^iu^j \, ,
\eeq 
for all   $u=\sum_{i=1}^m u^i\p_{x_i} \in T_pM$, where
\beq
\psi_\lambda &=&\frac{\al_1[(\be_1+\be_3)\phi_1-\be_2\phi_2]
-\al_2(\al_1\be_2-\al_2\be_1)}{\al \phi}\\ \nb
 \psi_\theta &=&\frac{-\al_2[(\be_1+\be_3)\phi_1-\be_2\phi_2]
 +(\al_1+\al_3)(\al_1\be_2-\al_2\be_1)}{\al \phi}\\ \nb
  \psi_\omega &=&
\frac{(\al_1+\al_3)[\be_1(\phi_1+\phi_3)-\be_2\phi_2]
+\al_2[\al_2(\be_1+\be_3)-\be_2(\al_1+\al_3)]}{\al \phi}\, . \nb
  \eeq
\epro
\underline{{\it Proof}}\\ 
\vspace{0.15cm}

The product of the matrix-value functions $G$ and $G^{-1}$  block per block
 gives:
\beq
\left(
\begin{array}{ccc}
M_1 +M_3 & {}& M_2 \\
         &{} &    \\
	 M_2      & {}& M_1
	 \end{array}
	 \right)
\left(
\begin{array}{ccc} \Lambda & {} & \Theta \\
  & {} &{} \\
  \Theta & {} & \Omega
  \end{array}
  \right)
=\left(
\begin{array}{ccc}
(M_1 +M_3)\Lambda +M_2\Theta & {}& (M_1+M_3)\Theta
+M_2\Omega \\
         &{} &    \\
M_2\Lambda + M_1\Theta    & {}& M_2\Theta+M_1\Omega
\end{array}
\right)
\eeq
and so we have  the identities:
\beq
(M_1 +M_3)\Lambda +M_2\Theta &=& Id \label{e1}\\
(M_1+M_3)\Theta+M_2\Omega &=& 0   \label{e2}\\
M_2\Lambda + M_1\Theta  &=& 0    \label{e3} \\
 M_2\Theta+M_1\Omega &=& Id\ .   \label{e4}
\eeq
Furthermore, for any $ u\in T_pM$, since $(U;\, x_i,i=1,\cdots, m)$
is a normal coordinates system centred at $p$, we have  
$(M_1+M_3)(p,u)= \mu(\al_1+\al_3\, ,\, \be_1+\be_3\, ,\, u)\, ,\\
\quad M_2(p,u) = \mu(\al_2\, ,\, \be_2,\, u)\, ,
\quad M_1(p,u) = \mu(\al_1\, ,\, \be_1, u)\, $ ; where 
$u\equiv (u^i)_{i=1,\cdots,m}$~.
Then according to the system (\ref{Lp4}) and Lemma \ref{L2}, the
matrix-value functions $M_1$ and $(M_1+M_3)$ at $(p,u)$ are 
invertible. It follows that at $(p,u)$, the identities  
(\ref{e3}) and  (\ref{e2}) give respectively 
\beq
\Theta=-M_1^{-1}M_2\Lambda \label{e5}
\eeq 
and
\beq
\Theta=-(M_1+M_3)^{-1}M_2\Omega\label{e6}.
\eeq
 Combining the identities (\ref{e5}) and (\ref{e1}), we obtain
$$ (M_1+M_3 - M_2M_1^{-1}M_2)\Lambda=Id\, .$$
So $\Lambda(p,u)$ is invertible with
\beq\label{e8}
\Lambda(p,u) = (M_1+M_3 - M_2M_1^{-1}M_2)^{-1}_{|_{(p,u)}}\, .
\eeq
Next  we compute the elements of the matrix-value function \\
$(M_1+M_3 - M_2M_1^{-1}M_2)$ at $(p,u)$,
and we obtain
\beq\label{e12}
[(M_1+M_3)-M_2M_1^{-1}M_2]_{ij}&=& \lambda_1\delta_{ij}+\lambda_2u^iu^j
\eeq
where
\beq 
\lambda_1&=&\frac{\al}{\al_1}\quad \mbox{and} \\ \nb
\lambda_2&=&\frac{\phi_1[\al_1(\be_1+\be_3)-\al_2\be_2
-\phi_2\be_2] +\be_1\phi_2^2}{\al_1\phi_1}\, ,
\eeq
with
  \beq\label{e14}
  \lambda_1\neq 0\,  \mbox{ and }\,
   (\lambda_1+|u|^2\lambda_2)=\frac{\phi}{\phi_1}\neq 0\, .
 \eeq
   So by Lemma \ref{L2}, we obtain the inverse 
   $\Lambda=(\lambda^{ij})_{1\leq i,j\leq m}$ of\\ 
   $[(M_1+M_3)-M_2M_1^{-1}M_2]$ at $(p,u)$, with
   \beq\label{e13}
   \lambda^{ij}(p,u)&=&\frac{\delta_{ij}}{\lambda_1}
   -\frac{\lambda_2}
   {\lambda_1(\lambda_1+|u|^2\lambda_2)}u^iu^j
             \\ \nb
        &=& \frac{\al_1}{\al}\delta_{ij}
	  -\psi_\lambda u^iu^j\, .
    \eeq

Next, according to (\ref{e5}), we compute
\beq\label{e9}
\theta^{ij}(p,u)=-[M_1^{-1}M_2\Lambda]_{{ij}_{|_{(p,u)}}}\, ,
\eeq
and we obtain (\ref{2e}).

Furthermore by combining  (\ref{e6}) and (\ref{e4}) we obtain
\beq
[-M_2(M_1+M_3)^{-1}M_2+M_1]\Omega &=& Id\, .
\eeq
This  shows that the matrix-value function  
$\left[-M_2(M_1+M_3)^{-1}M_2+M_1\right]$ \\
is  invertible and
\beq\label{e7}
\Omega= [M_1-M_2(M_1+M_3)^{-1}M_2]^{-1}\, \mbox{ at } (p,u) .
\eeq
Finally, as in the proof of (\ref{e13}), we obtain 
\beq\label{e16}
{[M_1-M_2(M_1+M_3)^{-1}M_2]_{ij}}_{|_(p,u)}=\omega_1\delta_{ij}+\omega_2u^iu^j\, ,
\eeq
where  $\omega_1=\frac{\al}{\al_1+\al_3}\neq 0\, $   and
$\, \omega_2=\frac{(\phi_1+\phi_3)[\be_1(\al_1+\al_3)-\al_2\be_2-\be_2\phi_2]
                +\phi_2^2(\be_1+\be_3)}{(\al_1+\al_3)(\phi_1+\phi_3)}\quad$
		with
\beq\label{e15}
\omega_1+|u|^2\omega_2=\frac{\phi}{\phi_1+\phi_3}\neq 0\, .
\eeq
So by using again Lemma \ref{L2}, we prove (\ref{3e}).
                           
			                   $\hfill{\square}$

\brmq
The functions $\psi_\lambda$, $\psi_\theta$ and $\psi_\omega$ in
Proposition \ref{P4} only depend on the norms of the vectors 
$u \in T_pM$, since the same holds for  the  functions 
$\al_i\, \be_i\; ; i=1,2,3$.
\ermq
Besides we have the following lemma: 
\blem\label{pru1}
If $\al(t)\phi(t)\neq 0,\, \forall\, t\in \mbb{R}^+ $, then the
functions $\psi_\lambda$, $\psi_\theta$, $ \psi_\omega$ defined
respectively  in (\ref{1e}), (\ref{2e}) and (\ref{3e})
satisfy on $\mbb{R}^+$ the following identities:
\beq
\phi_2\psi_\lambda +\phi_1\psi_\theta & =& \frac{\al_1\be_2-\al_2\be_1}{\al}
\, ,\label{rr4} \\
(\phi_1+\phi_3)\psi_\lambda +\phi_2\psi_\theta
               &=& \frac{\al_1(\be_1+\be_3)-\al_2\be_2}{\al}\, ,\label{rr5} \\
\phi_2\psi_\theta +\phi_1\psi_\omega & =& 
\frac{(\al_1+\al_3)\be_1-\al_2\be_2}{\al}\, ,\label{rr6} \\
(\phi_1+\phi_3)\psi_\theta +\phi_2\psi_\omega & =& 
\frac{(\al_1+\al_3)\be_2-\al_2(\be_1+\be_3)}{\al}\, . \label{rr7}
\eeq
\elem

The proof of the identities of Lemma \ref{pru1} is not very 
difficult and can be obtained by straightforward computations  .
\vspace{0.15 cm}

\bpro\label{P5}
If $G$ is nondegenerate,   the elements of the 
matrix-value functions in (\ref{e})  are defined by
\beq
\lambda^{ij}(x,u) = \frac{\alpha_1}{\alpha} g^{ij} 
-\psi_\lambda u^iu^j
\eeq
\beq
\theta^{ij}(x,u) =-\frac{\alpha_2}{\alpha} g^{ij} 
-\psi_\theta u^iu^j
\eeq
\beq
\omega^{ij}(x,u) =\frac{\alpha_1+\alpha_3}{\alpha} g^{ij} 
-\psi_\omega u^iu^j\, ;
\eeq
for any   
$(x\, ,u)\in\pi^{-1}( U)$, with $u=
\sum_{i=1}^m u^i\p_{x_i} \in T_xM$;
where 
$\left(g^{ij}\right)_{1\leq i,j\leq m}$ denotes the inverse of 
$g\equiv \left(g_{ij}\right)_{1\leq i,j\leq m}$ with 
$g_{ij}= g(\p_{x_i}, \p_{x_j})$.
\epro
\underline{\it Proof}\\
Let us set
$$L=
\left(
\begin{array}{ccc}
(M_1 +M_3) & {}& M_2 \\
         &{} &    \\
	 M_2      & {}& M_1
	 \end{array}
	 \right)\, 
 \left(
\begin{array}{ccc}
\left(\frac{\alpha_1}{\alpha} g^{ij} -\psi_\lambda u^iu^j\right)_{1\leq i,j\leq m} & {}&
\left(-\frac{\alpha_2}{\alpha} g^{ij} -\psi_\theta u^iu^j\right)_{1\leq i,j\leq m}
\\
         &{} &    \\
\left(-\frac{\alpha_2}{\alpha} g^{ij} -\psi_\theta u^iu^j\right)_{1\leq i,j\leq m} & {}&
\left(\frac{\alpha_1+\alpha_3}{\alpha} g^{ij}
	        -\psi_\omega u^iu^j\,\right)_{1\leq i,j\leq m}
	 \end{array}
	 \right)\, ,
$$
with $L=(L_{ij})_{1\leq i,j\leq 2m}$.

It suffices to show that $L_{ij}=\delta_{ij}\; ; \mbox{ for }i,j=1,\dots, 2m$.
Actually, we have  for $i\, ,\, j=1,\cdots m$ :
\beq\label{i1}
L_{ i j}&=&\sum_{k=1}^m[(\al_1+\al_3)g_{ik}+(\be_1+\be_3)g(\p_{x_i},u)
g(\p_{x_k},u)]
              [\frac{\al_1}{\al}g^{kj}-\psi_\lambda u^ku^j]\nb \\ 
	      &&+\sum_{k=1}^m[\al_2g_{ik}+\be_2g(\p_{x_i},u)g(\p_{x_k},u)]
	      [-\frac{\al_2}{\al}g^{kj}-\psi_\theta u^ku^j] \\ \nb
&=& \frac{\al_1(\al_1+\al_3)}{\al}\sum_{k=1}^mg_{ik}g^{kj}
-(\al_1+\al_3)\psi_\lambda u^j\sum_{k=1}^mg_{ik}u^k\\ \nb
&&+\frac{\al_1(\be_1+\be_3)}{\al}g(\p_{x_i},u)\sum_{k=1}^mg(\p_{x_k},u)g^{kj}
-(\be_1+\be_3)\psi_\lambda g(\p_{x_i},u)u^j\sum_{k=1}^mg(\p_{x_k},u)u^k\\ \nb
&&-\frac{\al_2^2}{\al}\sum_{k=1}^mg_{ik}g^{kj}-\al_2\psi_\theta u^j
\sum_{k=1}^mg_{ik}u^k\\
                       \nb
&& -\frac{\al_2\be_2}{\al}g(\p_{x_i},u)\sum_{k=1}^mg(\p_{x_k},u)g^{kj}
-\be_2\psi_\theta g(\p_{x_i},u)u^j\sum_{k=1}^m g(\p_{x_k},u)u^k \\ \nb
&=& \frac{\al_1(\al_1+\al_3)}{\al}\delta_{ij}
-(\al_1+\al_3)\psi_\lambda u^j g(\p_{x_i},u)\\ \nb
&&+\frac{\al_1(\be_1+\be_3)}{\al}g(\p_{x_i},u)u^j
-(\be_1+\be_3)\psi_\lambda g(\p_{x_i},u)u^jg(u,u)\\ \nb
&&-\frac{\al_2^2}{\al}\delta_{ij}-\al_2\psi_\theta u^j
g(\p_{x_i},u)\\
                       \nb
 && -\frac{\al_2\be_2}{\al}g(\p_{x_i},u)u^j
-\be_2\psi_\theta g(\p_{x_i},u)u^j g(u,u) \\ \nb
&=& \delta_{ij}+[\frac{\al_1(\be_1+\be_3)-\al_2\be_2}{\al}
  -(\phi_1+\phi_3)\psi_\lambda-\phi_2\psi_\theta]
   g(\p_{x_i},u)u^j\\ \nb
L_{ i j}&=&\delta_{ij}\; \mbox{ by } (\ref{rr5})\; ,
 \eeq
 \beq\label{i2}
 L_{\{i+m\} j}&=&\sum_{k=1}^m[\al_2g_{ik}+\be_2g(\p_{x_i},u)g(\p_{x_k},u)]
               [\frac{\al_1}{\al}g^{kj}-\psi_\lambda u^ku^j]\\ \nb
	                     &&+\sum_{k=1}^m[\al_1g_{ik}
			     +\be_1g(\p_{x_i},u)g(\p_{x_k},u)]
			 [-\frac{\al_2}{\al}g^{kj}-\psi_\theta u^ku^j] \\ \nb
&=& (\frac{\al_1\be_2-\al_2\be_1}{\al}-\phi_2\psi_\lambda-\phi_1\psi_\theta)
g(\p_{x_i},u)u^j
    \\ \nb
L_{\{i+m\} j}&=& 0\; \mbox{ by } (\ref{rr4})\; ,
\eeq
\newpage
\beq\label{i3}
L_{\{i+m\}\{j+m\}}&=&\sum_{k=1}^m[\al_2g_{ik}+\be_2g(\p_{x_i},u)
           g(\p_{x_k},u)]
              [-\frac{\al_2}{\al}g^{kj}-\psi_\theta u^ku^j]\\ \nb
	 &&+\sum_{k=1}^m[\al_1g_{ik}+\be_1g(\p_{x_i},u)g(\p_{x_k},u)]
	 [\frac{(\al_1+\al_3)}{\al}g^{kj}-\psi_\omega u^ku^j] \\
	 \nb
	    &=& \delta_{ij}
+[\frac{(\al_1+\al_3)\be_1-\al_2\be_2}{\al}
-(\phi_1\psi_\omega+\phi_2\psi_\theta)]g(\p_{x_i},u)u^j\\ \nb
L_{\{i+m\}\{j+m\}} &=& \delta_{ij}\; \mbox{ by } (\ref{rr6})\; ,
 \eeq
\beq\label{i4}
L_{ i\{j+m\}}&=&\sum_{k=1}^m[(\al_1+\al_3)g_{ik}
+(\be_1+\be_3)g(\p_{x_i},u)g(\p_{x_k},u)]
              [-\frac{\al_2}{\al}g^{kj}-\psi_\theta u^ku^j]\nb\\ 
&&+\sum_{k=1}^m[\al_2g_{ik}+\be_2g(\p_{x_i},u)g(\p_{x_k},u)]
[\frac{(\al_1+\al_3)}{\al}g^{kj}-\psi_\omega u^ku^j] \\ \nb
          &=& 
[\frac{(\al_1+\al_3)\be_2-\al_2(\be_1+\be_3)}{\al}
 -\phi_2\psi_\omega-(\phi_1+\phi_3)\psi_\theta]g(\p_{x_i},u)u^j\\ \nb
L_{ i\{j+m\}} &=& 0\; \mbox{ by } (\ref{rr7})\; .
\eeq
Hence $L_{ij} =\delta_{ij}\; \mbox{ for } i,j=1,\cdots, 2m\, $;
as stated.

\hfill{$\square$}
 \section{Levi-Civita connection of a nondegene\-rate $g$-natural metric}
\label{S4}

In \cite{AS1}, the authors have given explicitly (with some sign and 
parenthesis misprints) the Levi-Civita connection in the case of  Riemannian 
$g$-natural me\-trics. In the  following we determine the Levi-Civita connection
for a nondege\-nerate $g$-natural metric in general by using the inverse
formula  of nondege\-nerate $g$-natural metrics.
\bnota\label{not}
For a Riemannian manifold $(M,g)$,  we set :
\beq
\begin{array}{lcl}
T^1(u;X_x,Y_x)=R(X_x,u)Y_x, & & T^2(u;X_x,Y)=R(Y_x,u)X_x\, ,\\
T^3(u;X_x,Y_x)=R(X_x,Y_x)u, & & T^4(u;X_x,Y_x)=g(R(X_x,u)Y_x,u)u\, ,\\
T^5(u;X_x,Y_x)=g(X_x,u)Y_x, & & T^6(u;X_x,Y_x)=g(Y_x,u)X_x\, , \\
T^7(u;X_x,Y_x)=g(X_x,Y_x)u, & & T^8(u;X_x,Y_x)=g(X_x,u)g(Y_x,u)u\, ,
\end{array}
\eeq
where $(x,u)\in TM\, $,  $\, X_x,Y_x\in T_xM\, $ and $\, R\, $ is the Riemannian
curvature of $g\, $.
\enota
Let  $ \n$ be the   Levi-Civita connection  of $g$ and $\bar{\n}$
the  Levi-Civita connection of a nondegenerate  $g$-natural metric $G$ 
defined by the functions $\al_i,\ \be_i, \ i=1,2,3$ in 
Proposition \ref{P1}. We have:

\bpro\label{P6}
Let $(x,u)\in TM$ and $X,Y\in \mf{X}(M)$, we have
\beq
\left(\bar{\n}_{X^h}Y^h\right)_{(x,u)} &=&\left(\n_XY\right)_{(x,u)}^h
              + h\{A(u;X_x,Y_x)\}+v\{B(u;X_x,Y_x)\}\quad   \label{h1}\\ 
\left(\bar{\n}_{X^h}Y^v\right)_{(x,u)} &=&\left(\n_XY\right)_{(x,u)}^v 
              + h\{C(u;X_x,Y_x)\}+v\{D(u;X_x,Y_x)\}\quad \label{h2} \\ 
\left(\bar{\n}_{X^v}Y^h\right)_{(x,u)} &=& h\{C(u;Y_x,X_x)\}+v\{D(u;Y_x,X_x)\}
\label{h3}\\
\left(\bar{\n}_{X^v}Y^v\right)_{(x,u)} &=& h\{E(u;Y_x,X_x)\}+v\{F(u;Y_x,X_x)\}
\label{h4}
\eeq
where
$P(u;X_x,Y_x)=\sum_{i=1}^8f^P_i(|u|^2)T^i(u;X_x,Y_x),\quad \mbox{for }
P=A,B,C,D,E,F$;
with
\beq
\begin{array}{lcl}\label{A}
f_1^A=\, f_2^A =-\frac{\al_1\al_2}{2\al}\, , & & f_3^A = 0\, ,\\
        & &    \\
f_4^A=\al_2\psi_\lambda\, ,& & f_5^A=f_6^A = \frac{\al_2(\be_1+\be_3)}{2\al}\, ,\\
        & &    \\
 f_7^A= (\al_1+\al_3)^\pr\frac{\phi_2}{\phi}\, ,
   && f_8^A = (\be_1+\be_3)^\pr\frac{\phi_2}{\phi}
           +(\be_1+\be_3)\psi_\theta\,  ;
 \end{array}
 \eeq
\beq
\begin{array}{lcl}\label{B}
f_1^B=\frac{\al_2^2}{\al}\, ,& & f_2^B = 0\, ,\\
        & &    \\
f_3^B=-\frac{\al_1(\al_1+\al_3)}{2\al}\, ,& & f_4^B=\al_2\psi_\theta\, ,\\
		        & &    \\
 f_5^B=f_6^B = -\frac{(\al_1+\al_3)(\be_1+\be_3)}{2\al}\, ,
        & & f_7^B= -(\al_1+\al_3)^\pr\frac{(\phi_1+\phi_3)}{\phi}\, ,\\
                          & &    \\
	 f_8^B =-(\be_1+\be_3)^\pr \frac{(\phi_1+\phi_3)}{\phi}
         +(\be_1+\be_3)\psi_\omega\, ; & & 
         \end{array}
	 \eeq

\beq
\begin{array}{lcl}\label{C}
f_1^C=0 \, ,& & f_2^C = -\frac{\al_1^2}{2\al}\, ,\\
        & &    \\
	f_3^C=0\, ,& & f_4^C=\frac{\al_1\psi_\lambda}{2}\, ,\\
	& &    \\
	f_5^C = +\frac{\al_1(\be_1+\be_3)}{2\al}\, ,
	    & & f_6^C= (\al_1+\al_3)'\frac{\al_1}{\al}
	           -\frac{\al_2}{2\al}(2\al_2'-\be_2)\, ,\\
		          & &    \\
f_7^C = \frac{(\be_1+\be_3)\phi_1}{2\phi}
        +\frac{1}{2}(2\al_2^\pr-\be_2)\frac{\phi_2}{\phi}\, ,
 & & f_8^C= (\be_1+\be_3)^\pr\frac{\phi_1}{\phi}
      -\psi_\lambda[(\al_1+\al_3)^\pr+\frac{(\be_1+\be_3)}{2}]\\ 
                  & &\hspace{1cm}
               -\frac{1}{2}(2\al_2^\pr-\be_2)\psi_\theta\, ;
	\end{array}
\eeq
\beq
\begin{array}{lcl}\label{D}
f_1^D = 0\, , & & f_2^D = \frac{\al_1\al_2}{2\al}\, ,\\
       & &  \\
       f_3^D = 0\, ,& & f_4^D =\frac{\al_1}{2}\psi_\theta\, ,\\
       & &  \\
       f_5^D =- \frac{\al_2(\be_1+\be_3)}{2\al}\, ,
            & & f_6^D= -(\al_1+\al_3)'\frac{\al_2}{\al}
	   +\frac{(2\al_2^\pr-\be_2)(\al_1+\al_3)}{2\al}\, , \\
	      & &  \\
 f_7^D= -\frac{(\be_1+\be_3)\phi_2}{2\phi}
          -\frac{1}{2}(2\al_2^\pr-\be_2)\frac{(\phi_1+\phi_3)}{\phi}\, ,
	   & & f_8^D=-(\be_1+\be_3)^\pr\frac{\phi_2}{\phi}
	        -[(\al_1+\al_3)^\pr+\frac{\be_1+\be_3}{2}]\psi_\theta\\
	 &&   \hspace{1cm}  -\frac{1}{2}(2\al_2^\pr-\be_2)\psi_\omega\, ;
	 \end{array}
 \eeq
\beq
\begin{array}{lcl}\label{E}
f_1^E=f_2^E=f_3^E=f_4^E=0 \, ,& & 
f_5^E=f_6^E=(\al_2^\pr+\frac{1}{2}\be_2)\frac{\al_1}{\al}-\al_1^\pr\frac{\al_2}{\al}
             \, ,   \\ 
		& &  \\ 
f_7^E =\be_2\frac{\phi_1}{\phi}
     -(\be_1-\al_1^\pr)\frac{\phi_2}{\phi} \, ,& &f_8^E = 
     2\be_2^\pr\frac{\phi_1}{\phi}
 -\be_1^\pr\frac{\phi_2}{\phi}-(2\al_2^\pr+\be_2)\psi_\lambda
     -2\al_1^\pr\psi_\theta\, ;
     \end{array}
\eeq
\beq
\begin{array}{lcl}\label{F}
f_1^F=f_2^F=f_3^F=f_4^F=0\, , & &
f_5^F=f_6^F=-(\al_2^\pr+\frac{1}{2}\be_2)\frac{\al_2}{\al}
+\al_1^\pr\frac{(\al_1+\al_3)}{\al}
            \, ,    \\ 
		                & &  \\ 
f_7^F =(\be_1-\al_1')\frac{(\phi_1+\phi_3)}{\phi} 
-\frac{\be_2\phi_2}{\phi} \, ,& &
f_8^F = \be_1'\frac{(\phi_1+\phi_3)}{\phi}
    -2\be_2^\pr \frac{\phi_2}{\phi}-(2\al_2^\pr+\be_2)\psi_\theta
		-2\al_1^\pr\psi_\omega\; .
      \end{array}
\eeq
\epro
\underline{{\it Proof}}
\vspace{0.25cm}

We  prove only (\ref{h3}), the proof of the other being the same.  
Let us set 
\beq
X= \sum_{i=1}^mX^i\p_{x_i}\, , \quad
Y= \sum_{i=1}^mY^i\p_{x_i} \, , \quad
u= \sum_{i=1}^mu^i\p_{x_i} \, , \quad
\eeq
\beq
\bar{\n}_{X^v}Y^h&=&\sum_{i=1}^md_i\p_{x_i}^h+\sum_{i=1}^md_{m+i}\p_{x_i}^v
   \label{naa}\\                  
s_i&=&G\left(\bar{\n}_{X^v}Y^h,\p_{x_i}^h\right)\quad \mbox{and}\\ 
s_{m+i}&=&G\left(\bar{\n}_{X^v}Y^h,\p_{x_i}^v\right) \, . 
\eeq
 Koszul's formula  gives
\beq
s_i&=&\frac{1}{2}
\left\{X^v.G\left(Y^h,\p_{x_i}^h\right)+Y^h.G\left(\p_{x_i}^h,X^v\right)
-\p_{x_i}^h.G\left(X^v,Y^h\right)\right.\\ \nb
 {} &{}& \left.  +G\left(\p_{x_i}^h,\left[X^v,Y^h\right]\right)
 -G\left(Y^h,\left[X^v,\p_{x_i}^h\right]\right)
 -G\left(X^v,\left[Y^h,\p_{x_i}^h\right]\right)
 \right\},\nb 
 \eeq
then by using  Proposition \ref{P1}, Lemma \ref{L1}  and Lemma \ref{L1p},
we obtain
\beq
s_i&=& (\al_1+\al_3)'g(X,u)g(Y,\p_{x_i})
      +(\be_1+\be_3)'g(X,u)g(Y,u)g(\p_{x_i},u)\quad \\  \nb
        &{}&+\frac{\be_1+\be_3}{2}g(X,Y)g(\p_{x_i},u)
	       +\frac{\be_1+\be_3}{2}g(Y,u)g(X,\p_{x_i})\\ \nb
	         &{}&     +\frac{\al_1}{2}g(R(Y,\p_{x_i})u,X)
		 \nb
\eeq
and similarly
\beq
s_{m+i}=\frac{1}{2}(2\al_2'-\be_2)g(X,u)g(Y,\p_{x_i})
        -\frac{1}{2}(2\al_2'-\be_2)g(X,Y)g(u,\p_{x_i})
\eeq
By setting $\mathbf{d}=(d_i)_{1\leq i\leq 2m}$ and 
$\mathbf{s}=(s_i)_{1\leq i\leq 2m}$, we have
$\mathbf{ d}=G^{-1}\mathbf{s}$
 (Matrix-value function of  $G^{-1}$ with the column vector 
 $\mathbf{s}$ as argument).

Then by using the expression of $G^{-1}$ in  Proposition \ref{P5}, we obtain 
\beq\label{lbi}
d_i&=& \frac{\al_1^2}{2\al}\{R(u,X)Y\}^i
  -\frac{\al_1\psi_\lambda}{2}g(R(Y,u)u,X)u^i\ \\ \nb
  &&+[(\al_1+\al_3)'\frac{\al_1}{\al}
  -\frac{\al_2}{2\al}(2\al_2'-\be_2)]g(X,u)Y^i
      \\ \nb
   &&+\frac{\al_1(\be_1+\be_3)}{2\al}g(Y,u)X^i  \\ \nb
   &&+[\frac{1}{2}(\be_1+\be_3)\frac{\phi_1}{\phi}
   +\frac{1}{2}(2\al_2'-\be_2)\frac{\phi_2}{\phi}]
      g(X,Y)u^i \\ \nb
     &&+\{(\be_1+\be_3)'\frac{\phi_1}{\phi}
     -[(\al_1+\al_3)'+\frac{\be_1+\be_3}{2}]\psi_\lambda\\ \nb
  &&  -\frac{1}{2}\psi_\theta(2\al_2'-\be_2)\}g(X,u)g(Y,u)u^i \nb
    \eeq
and
\beq\label{lbni}
d_{m+i}&=& \frac{\al_1\al_2}{2\al}\{R(X,u)Y\}^i
 +\frac{\al_1\psi_\theta}{2}g(R(X,u)Y,u)u^i\\ \nb
  &&+[-(\al_1+\al_3)'\frac{\al_2}{\al}
   +\frac{(2\al_2'-\be_2)(\al_1+\al_3)}{2\al}]g(X,u)Y^i\\ \nb
    &&-\al_2\frac{(\be_1+\be_3)}{2\al}g(Y,u)X^i \\ \nb
     &&+[-\frac{(\be_1+\be_3)\phi_2}{2\phi}
      -\frac{1}{2}(2\al_2'-\be_2)\frac{\phi_1+\phi_3}{\phi}]g(X,Y)u^i\\ \nb
       &&+\{-(\be_1+\be_3)'\frac{\phi_2}{\phi}
        -[(\al_1+\al_3)'+\frac{(\be_1+\be_3)}{2}]\psi_\theta\\ \nb
	 && -\frac{1}{2}(2\al_2'-\be_2)\psi_\omega\}g(X,u)(Y,u)u^i\, ,\nb
	  \eeq
where for all $W\ \in\ \mathfrak{X}(M)\ $, $\{W\}^i$ are the components of 
 $W$ in the coordinates system $(U;x_{i},\, i=1,\cdots,m) $. So 
 according  to (\ref{naa}), 
 the proof of  (\ref{h3}) is completed. 

$\hfill{\square}$

\section{$g$-Natural metrics  with  constant sectional curvature}\label{S5}
\subsection{Riemannian curvature of nondegenerate $g$-natural
metrics}
\subsubsection*{Some notations and properties of $F$-tensor fields }
Fix $(x,u)\in TM$ and a system of normal coordinates \\
$S:= (U\, ;\, x_i\, ,i=1,\cdots, m)$ of $(M,g)$ centred at $x$. Then we can 
define on $U$ the vector field $\mathbf{U} := \sum_iu^i\frac{\p}{\p x_i}$,
where $(u^1,\cdots,u^m)$ are the coordinates of $u\in T_xM$ with respect
to its basis $ (\left(\frac{\p}{\p x_i}\right)_x; \, i=1,\cdots,m )$.

Let $P$ be an $F$-tensor field of type $(r,s)$ on $M$. Then, on $U$, we can 
define an $(r,s)$-tensor field $P_u^S$ (or $P_u$ if there is no risk of 
confusion), associated to $u$ and $S$, by
\beq
P_u(X_1,\cdots,X_s) := P(\mathbf{U}_z; X_1,\cdots,X_s)\, ,
\eeq
for all $(X_1,\cdots,X_s)\in T_zM,\; \forall z\in U$.

On the other hand, if we fix $x\in M$ and $s$ vectors $X_1,\cdots,X_s$ in 
$T_xM$, then we can define a $C^\infty$-mapping
$P_{(X_1,\cdots,X_s)}: T_xM\rightarrow \otimes^rT_xM$, associated to 
$(X_1,\cdots,X_s)$, by 
\beq
P_{(X_1,\cdots,X_s)}(u) :=P(u;\, X_1,\cdots,X_s)\, ,
\eeq
for all $u\in T_xM$.

Let $s>t$ be two non-negative integers, $T$ be a $(1,s)$-tensor field on $M$
and $P^T$ be an $F$-tensor field, of type $(1,t)$, of the form
\beq
P^T(u; X_1,\cdots, X_t) = T(X_1,\cdots, u,\cdots, u,\cdots, X_t) ,
\eeq
for all $(u; X_1,\cdots, X_t)\in TM\oplus\cdots\oplus TM $, i.e., $u$ appears
$s-t$ times at positions $i_1,\cdots,i_{s-t}$ in the expression of $T$. Then
\begin{itemize}
\item[{\bf -}] $P_u^T$ is a $(1,t)$-tensor field on a neighborhood $U$ of
$x$ in $M$,\\ for all $u\in T_xM$~;
\item[{\bf -}] $P_{(X_1,\cdots,X_t)}^T$ is a $C^\infty$-mapping 
$T_xM\rightarrow T_xM$, for all $X_1,\cdots,X_t$ in $T_xM$. 
\end{itemize}
Furthermore, we have 
\blem\cite{AS2}
\begin{enumerate}
\item[1)] The covariant derivative of $P_u^T$, with respect to the Levi-Civita
connection of $(M,g)$, is given by :
\beq
\left(\n_X P_u^T\right)(X_1,\cdots,X_t) = (\n_X T)
(X_1,\cdots,u,\cdots,u, X_t),
\eeq
for all vectors $X,X_1,\cdots,X_t$ in $T_xM$, where $u$ appears at
positions \\ $i_1,\cdots,i_{s-t}$ in the right-hand side of the 
preceding formula.
\item[2)] The differential of $P_{(X_1,\cdots,X_t)}^T$, at $u\in T_xM$, is
given by :
\beq
d\left(P_{(X_1,\cdots,X_t)}^T\right)_u(X) &=&
T(X_1,\cdots,X,\cdots, u,\cdots,X_t) +\cdots\\ \nb 
 & &   + T(X_1,\cdots,u,\cdots, X,\cdots,X_t),
\eeq
for all $X\in T_xM$.
\end{enumerate}
\elem
Furthermore, in \cite{AS2} the authors gave the expressions determining the\\
Riemannian curvature $\bar{R}$ of any Riemannian $g$-natural metric $G$ on
$TM$ (up to a misprint in the vertical component of the expression of 
$\bar{R}\left(X^h,Y^h\right)Z^h$, in which
$\left(\n_YA_u\right)(X,Z)$ should
be written $\left(\n_YB_u\right)(X,Z) )$. Their formulas remain the same
if we replace a Riemannian $g$-natural metric by a nondegene\-rate 
$g$-natural metric  on $TM$. Indeed, a similar proof as that in \cite{AS2}
gives~: 
\bpro
The Riemannian curvature $\bar{R}$ of a  nondegenerate  $g$-natural metric $G$ 
is completely defined by
\beq\label{rli}
\bar{R}\left(X^h,Y^h\right)Z^h&=& h\{R(X,Y)Z\}\\ \nb 
      &&+h\{(\n_XA_u)(Y,Z)-(\n_YA_u)(X,Z)\\ \nb
      &&+A(u;X,A(u;Y,Z))-A(u;Y,A(u;X,Z))\\ \nb
      &&+C(u;X,B(u;Y,Z))-C(u;Y,B(u;X,Z))\\ \nb
      &&+C(u;Z,R(X,Y)u)\}\\ \nb
      &&+v\{(\n_XB_u)(Y,Z)-(\n_YB_u)(X,Z)\\ \nb
      &&+B(u;X,A(u;Y,Z))-B(u;Y,A(u;X,Z))\\ \nb
      &&+D(u;X,B(u;Y,Z))-D(u;Y,B(u;X,Z))\\ \nb
      &&+ D(u;Z,R(X,Y)u)\}\, ,\nb
\eeq
\beq\label{rlii}
\bar{R}\left(X^h,Y^h\right)Z^v&=&h\{(\nabla_XC_u)
                    \left(Y,Z\right)-(\nabla_YC_u)\left(X,Z\right)
	\\ \nonumber
	                 &&+A\left(u;X,C\left(u;Y,Z\right)\right)
	                -A\left(u;Y,C\left(u;X,Z\right)\right)
	                +C\left(u;X,D\left(u;Y,Z\right)\right)\\
	                \nonumber
	                &&-C\left(u;Y,D\left(u;X,Z\right)\right)
	               +E\left(u;R\left(X,Y\right)u,Z\right)\}\\ \nb
	              && +v\{R(X,Y)Z
		       +\left(\nabla_XD_u\right)\left(Y,Z\right)
	                 -\left(\nabla_YD_u\right)\left(X,Z\right)\\ \nb
	               && +B\left(u;X,C\left(u;Y,Z\right)\right)
	              -B\left(u;Y,C\left(u;X,Z\right)\right)\\
                             \nonumber
                        &&+D\left(u;X,D\left(u;Y,Z\right)\right)
                       -D\left(u;Y,D\left(u;X,Z\right)\right)
                       +F(u;R\left(X,Y\right)u,Z)\}\, , \nonumber
  \eeq
\beq\label{rliii}
\bar{R}\left(X^h,Y^v\right)Z^h&=&h\{(\nabla_XC_u)\left(Z,Y\right)
                                 +A(u;X,C\left(u;Z,Y\right))\\ \nonumber
	                     && +C(u;X,D\left(u;Z,Y\right))
	                    -C(u;A\left(u;X,Z\right),Y) \\ \nonumber
	                    &&  -E(u;Y,B\left(u;X,Z\right))
	                    -d\left(A_{\left(X,Z\right)}\right)_u
	                          \left(Y\right)\}\\ \nonumber
	                   &&+v\{\left(\nabla_XD_u\right)\left(Z,Y\right)
	                         +B(u;X,C\left(u;Z,Y\right))
	                      +D(u;X,D\left(u;Z,Y\right))\\ \nonumber
                           &&-D(u;A\left(u;X,Z\right),Y)
			      -F(u;Y,B\left(u;X,Z\right))
			      -d\left(B_{(X,Z)}\right)_u(Y)\}\, ,
\eeq
\beq\label{rliv}
\bar{R}\left(X^h,Y^v\right)Z^v&=&h\{(\nabla_XE_u)\left(Y,Z\right)
                                 +A(u;X,E\left(u;Y,Z\right))\\ \nonumber
                            && +C(u;X,F\left(u;Y,Z\right))
                             -C(u;C\left(u;X,Z\right),Y)\\ \nonumber
                            &&  -E(u;Y,D\left(u;X,Z\right))
                            -d\left(C_{\left(X,Z\right)}\right)_u
                            \left(Y\right)\}\\ \nonumber
                           &&+ v\{\left(\nabla_XF_u\right)\left(Y,Z\right)
                            +B(u;X,E\left(u;Y,Z\right))
                            +D(u;X,F\left(u;Y,Z\right))\\ \nonumber
                            &&-D(u;C\left(u;X,Z\right),Y)
                             -F(u;Y,D\left(u;X,Z\right))
                             -d\left(D_{\left(X,Z\right)}\right)_u
                             \left(Y\right)\}\, ,\\ \nonumber
 \eeq
\beq\label{rlv}
\bar{R}\left(X^v,Y^v\right)Z^h&=&h\{d\left(C_{\left(Z,Y\right)}\right)_u
                                 \left(X\right)
                             -d\left(C_{\left(Z,X\right)}\right)_u
                              \left(Y\right)\\ \nonumber
                           && +C(u;C\left(u;Z,Y\right),X)
                           -C(u;C\left(u;Z,X\right),Y) \\ \nonumber
                          && +E(u;X,D\left(u;Z,Y\right))
                         -E(u;Y,D\left(u;Z,X\right))\}\\ \nonumber
                          &&+v\{d\left(D_{\left(Z,Y\right)}\right)_u
                            \left(X\right)
                         -d\left(D_{\left(Z,X\right)}\right)_u
                                          \left(Y\right)
		+D(u;C\left(u;Z,Y\right),X)\\ \nonumber
                         &&-D(u;C\left(u;Z,X\right),Y)
                         +F(u;X,D\left(u;Z,Y\right))
                        -F(u;Y,D\left(u;Z,X\right))\}\, ,
\eeq
\beq\label{rlvi}
\bar{R}\left(X^v,Y^v\right)Z^v&=&h\{d\left(E_{\left(Y,Z\right)}\right)_u
                                  \left(X\right)
                            -d\left(E_{\left(X,Z\right)}\right)_u
                                 \left(Y\right)
                              +C(u;E\left(u;Y,Z\right),X)\\ \nonumber
                               &&-C(u;E\left(u;X,Z\right),Y)
                              + E(u;X,F\left(u;Y,Z\right))
                           -E(u;Y,F\left(u;X,Z\right))\}\\ \nonumber
                           &&+v\{d\left(F_{\left(Y,Z\right)}\right)_u
                                      \left(X\right)
                              -d\left(F_{\left(X,Z\right)}\right)_u
                                     \left(Y\right)
                               +D(u;E\left(u;Y,Z\right),X)\\ \nonumber
                               &&-D(u;E\left(u;X,Z\right),Y)
                               +F(u;X,F\left(u;Y,Z\right))
                              -F(u;Y,F\left(u;X,Z\right))\}\, ,
\eeq
for all $x\in M$ and $X,Y,Z\in T_xM$, where the lifts are taken at $u\in T_xM$
and $R$ is the Riemannian curvature of $g$.
\epro
\vspace{0.15cm}

\brmq\label{R2}
Let $P=\sum_{i=5}^8f_i^PT^i\, , \quad Q=\sum_{i=5}^8f_i^QT^i\,$ 
be $F$-tensors.\\
For $(x,u)\in TM$ and $X\, ,\, Y\, , Z\, \in 
T_xM\,$, we have
\beq\label{p2}
P(u;\ X,Q(u;\ Y,Z))-P(u;\ Y,Q(u;\ X,Z)) &=&
\{a_1(P,Q)g(Y,Z)\\ \nb
&&+a_2(P,Q)g(Y,u)g(Z,u)\}X\\ \nb
&&-\{a_1(P,Q)g(X,Z)+a_2(P,Q)g(X,u)g(Z,u)\}Y\\ \nb
&&+a_3(P,Q)\{g(X,Z)g(Y,u)-g(Y,Z)g(X,u)\}u 
\eeq
 with
 \beq
 a_1(P,Q) &=& |u|^2f_6^Pf_7^Q\, ,\\
 a_2(P,Q)&=&f_6^P(f_6^Q+|u|^2f_8^Q)
    -(f_5^Pf_6^Q-f_6^Pf_5^Q)\, ,\\
 a_3(P,Q)&=&
   f_7^Pf_5^Q-(f_5^P+f_7^P+|u|^2f_8^P)f_7^Q\, ,
   \eeq
where $f_i^P\, , f_i^Q\, $ are differentiable functions on $\mbb{R}^+\ $ and
the $T^i$ are defined in Notation \ref{not} .
\ermq

In the sequel we shall consider only  Riemannian $g$-natural
metrics $G$ on $TM$.

\subsection{On the hereditary property of constant sectional curvature}
We prove the following result that improves \cite[theorem 0.3]{AS2}.
\vspace{0.15 cm}

\bpro\label{rp1}
 If $\, (TM,G)$ has constant sectional curvature then 
 $(M,g)$ is a flat Riemannian manifold. 
\epro
\underline{{ \it Proof}}
\vspace{0.25cm}

If  $\, (TM,G)$ has constant  sectional curvature $K$, then by 
\cite[theorem 0.3]{AS2}    $(M,g)$ has 
constant  sectional curvature  $k\in \mbb{R}$. Furthermore, since 
$(TM,G)$ has constant sectional curvature then its Riemannian
curvature $\bar{R}$ satisfies 
$\bar{R}(X^h,Y^h)Z^v_{|(x,u)}\in H_{(x,u)}TM$
for any $(x,u)\in TM, \, \mbox{ and }X,Y,Z\in \mf{X}(M)$.\\
Then by (\ref{rlii}), we have
\beq
R(X,Y)Z_{|x}&=& -[(\n_XD_u)(Y,Z)-(\n_YD_u)(X,Z)\\ \nb
     &&+B(u;X,C(u;Y,Z))-B(u;Y,C(u;X,Z)) \\ \nb
     &&+D(u;X,D(u;Y,Z))-D(u;Y,D(u;X,Z))\\ \nb
     && +F(u;R(X,Y)u,Z)]\, ,\\ \nb
     \forall (x,u)\in TM.\nb
\eeq
Thus $R(X,Y)Z_{|x}=0\, ,\, \forall x\in M$ (by taking $(x,u)=(x,0)\in TM$).\\
This means that  $k=0$.
                                                     
	     $\hfill{\square}$

In the following proposition, we investigated the $g$-natural metrics of\\
constant sectional curvature.
\vspace{0.15 cm}

\bpro\label{rp2}
For $dim\ M\geq 3$, the flat Riemannian  $g$-natural metrics are the only $g$-natural metrics on $TM$ that have constant sectional
curvature.
\epro
\underline{{ \it Proof}}
\vspace{0.25cm}

If $\, (TM,G)$ has constant sectional curvature  $K$, then
\beq\label{k00}
\bar{R}\left(X^h,Y^h\right)Z^h&=& K\left[G\left(Z^h,Y^h\right)X^h
             -G\left(X^h,Z^h\right)Y^h\right]\\ \nb
	     &=& K[(\al_1+\al_3)g(Z,Y)+
	          (\be_1+\be_3)g(Z,u)g(Y,u)]X^h\\ \nb
	&&-K[(\al_1+\al_3)g(X,Z)+
	        (\be_1+\be_3)g(X,u)g(Z,u)]Y^h. \nb
\eeq
So by  Proposition $\ref{rp1}$, we have $R\equiv 0$ and thus 
from the  formulas  (\ref{rli}) and  (\ref{p2}), we obtain  
\beq\label{k01}
\bar{R}\left(X^h,Y^h\right)Z^h&=&
       h\{A(u;X,A(u;Y,Z))-A(u;Y,A(u;X,Z))\\ \nb
      && +C(u;X,B(u;Y,Z))-C(u;Y,B(u;X,Z))\}\\ \nb
       &=&\{[a_1(A,A)+a_1(C,B)]g(Y,Z)
       +[a_2(A,A)+a_2(C,B)]g(Y,u)g(Z,u)\}X^h\\ \nb
    &&-\{[a_1(A,A)+a_1(C,B)]g(X,Z) \\ \nb
    &&      +[a_2(A,A)+a_2(C,B)]g(X,u)g(Z,u)\}Y^h\\ 
	   \nb
	&& +\{[a_3(A,A)+a_3(C,B)][g(X,Z)g(Y,u)
	      -g(Y,Z)g(X,u)]\}u^h \ .\nb
\eeq

Then,  let $(x,u)\in TM$ with $u\neq 0$ :
\begin{enumerate}
\item[1)] Since $\dim M \geq 3$,
there exists two vectors $X,\, Y\in T_xM$ such that the system
$(u,X,Y)$ is orthogonal. 
\newpage
So by (\ref{k00}) and (\ref{k01}), for  $Z=Y$, we obtain
respectively\\
$
\bar{R}\left(X^h,Y^h\right)Y^h= K(\al_1+\al_3)g(Y,Y)X^h
$
and \\
$
\bar{R}\left(X^h,Y^h\right)Y^h= [a_1(A,A)+a_1(C,B)]g(Y,Y)X^h\; ;
$ with $g(Y,Y)\neq 0$ and $X\neq 0$.  Hence  
\beq\label{k00a}
K(\al_1+\al_3)(t) =[a_1(A,A)+a_1(C,B)](t),\; \forall\, t> 0\, .
\eeq
\item[2)] Next, by choosing $Y=Z=u$ such as $u$ is orthogonal to a vector
$X\neq 0$
in $T_xM$, (\ref{k00}) gives 
\beq\label{k0b1}
\bar{R}\left(X^h,Y^h\right)Y^h= Kg(u,u)[(\al_1+\al_3)+g(u,u)(\be_1+\be_3)]X^h\; ,
\eeq
and  (\ref{k01}) gives
\beq\label{k0b2}
\bar{R}\left(X^h,Y^h\right)Y^h&=&g(u,u)[a_1(A,A)+a_1(C,B)\\
  & &       + g(u,u)(a_2(A,A)+a_2(C,B))] X^h .\nb
\eeq
Then, by (\ref{k0b1}) and (\ref{k0b2}), we have,
\beq
K[(\al_1+\al_3)+g(u,u)(\be_1+\be_3)] &=&
a_1(A,A)+a_1(C,B)\\ 
&& + g(u,u)[a_2(A,A)+a_2(C,B)] .\nb
\eeq

Thus, by (\ref{k00a}), we obtain
\beq\label{k00b}
[a_2(A,A)+a_2(C,B)](t) = K(\be_1+\be_3)(t)\; , \forall\; t>0\; .
\eeq
\item[3)] Furthermore, by choosing $Y=u$ and $X=Z\neq 0$ such as
$X$ and $u$ are orthogonal , (\ref{k00}) gives
$$\bar{R}\left(X^h,u^h\right)X^h = -K(\al_1+\al_3)g(X,X)u^h
$$
and (\ref{k01}) gives 
$$
\bar{R}\left(X^h,u^h\right)X^h =g(X,X)[ -(a_1(A,A)+a_1(C,B))+
                                g(u,u)(a_3(A,A)+a_3(C,B)]u^h\, .
$$
Then by (\ref{k00a}), we obtain 
\beq\label{k00c}
[a_3(A,A)+a_3(C,B)](t) =0\; , \forall\, t>0\, .
\eeq
And we deduce that the identities (\ref{k00a}), (\ref{k00b}) 
and (\ref{k00c}) are true for any $t\geq 0$, since the functions 
$\al_i\, ,\, \be_i\, ,i=1,2,3$ are smooth on $\mbb{R}^+$ .\\
Hence we have  
\beq\label{k02}
\left\{
\begin{array}{lcl}
a_1(A,A)+a_1(C,B) &=& K(\al_1+\al_3)\\
a_2(A,A)+a_2(C,B) &=& K(\be_1+\be_3)\\
a_3(A,A)+a_3(C,B) &=& 0\; .\\
\end{array}
\right.
\eeq
\end{enumerate}
But $(TM,G)$ is Riemannian i.e.; 
\beq
\left\{
\begin{array}{l}
\al_1 > 0\\
\al=\al_1(\al_1+\al_3)-\al_2^2 > 0
\end{array}
\right.\; ,
\eeq
 and then  
\beq
\left\{
\begin{array}{l}
\al_1 > 0\\
\al_1(\al_1+\al_3) > \al_2^2 
\end{array}
\right. \; ;
\eeq
so $ (\al_1+\al_3)>0\; $ .
Hence  according to the first equation of (\ref{k02}) which means that 
\beq
tf_6^A(t)f_7^A(t)+tf_6^C(t)f_7^B=K(\al_1+\al_3)(t),
\eeq
we obtain for $t=0$, $\ 0=K(\al_1+\al_3)(0)\ $, so $ K=0$.

                                   $\hfill{\square}$

If  $ (M,g)$ is  a flat Riemannian manifold and we choose
\beq
\left\{
\begin{array}{l}
\al_1\equiv 1,\\ 
\al_2=\al_3=\be_1=\be_2=\be_3\equiv 0
\end{array}
\right. ,
\mbox{ (Sasaki's metric) }
\eeq 
we obtain that
$(TM,G)$ is a flat Riemannian manifold . But it is not the only way to choose
the functions $\al_i,\ \be_i,\quad i=1,2,3$ for getting  $(TM,G)$ as a flat\\ 
Riemannian manifold. Actually we estabish a characterization of flat Riemannian 
$g$-natural metrics in what follows.
\subsection{Flat Riemannian $g$-natural metrics}
\blem\label{L5}
If $(TM,G)$ is a flat Riemannian manifold with $\dim M\geq 3$,
then
\begin{enumerate}
\item[a)] $\be_1+\be_3=0$,\label{l2i1}
\item[b)] $\al_1+\al_3 = constant > 0$,\label{l2i2}
\item[c)] $2\al_2'=\be_2$,     \label{l2i3}
\item[d)] $f_6^F=f_7^F=f_8^F=0$,\label{l2i4}
\end{enumerate}
where $\al_2'$ denote the first  derivative of $\al_2\, $.
\elem
\underline{{\it Proof}}:
\vspace{0.25cm}

If  $(TM,G)$ is flat Riemannian  then by \cite[page 36]{AS2},
we have $\be_1+\be_3=0$
and  $\al_1+\al_3=constant$. We have also  $\al_1+\al_3>0$ since  $\al>0$. 
Therefore we have the parts  $a)$ and $b)$ of Lemma \ref{L5}. Furthermore,\\
by \cite[Lemma 4.1]{AS3}, we have $\al_2^\pr-\be_2=0$, and
$A=B=C=D=0$.
\newpage

It remains to prove  $d)$.

Since $D=0$ then   (\ref{rlvi}) gives, 
\beq
\bar{R}_v^6(X,Y)Z &=&
\{[a_1(F,F)+f_7^F-f_6^F]g(Y,Z)\\ \nb
&& +[a_2(F,F)+f_8^F-2{f_6^F}^\pr]g(Y,u)g(Z,u)\}X\\ \nb 
&&-\{[a_1(F,F)+f_7^F-f_6^F]g(X,Z)\\ \nb 
&& +[a_2(F,F)+f_8^F-2{f_6^F}^\pr]g(X,u)g(Z,u)\}Y\nb\\
&&+[a_3(F,F)+2{f_7^F}^\pr-f_8^F]\{ g(X,Z)g(Y,u)-g(Y,Z)g(X,u)\}u.\nb
\eeq
where $\bar{R}_v^6(X,Y)Z$ is the vertical component of 
$\bar{R}\left(X^v,Y^v\right)Z^v$.
Since\\  $dim\ M \geq 3$  as before,
$\bar{R}_v^6(X,Y)Z=0\quad X,Y,Z\in T_xM,$ implies
\beq\label{S}
\left\{
\begin{array}{lcl}
tf_6f_7 +f_7-f_6 &=& 0 \\
f_6^2 +tf_6f_8 +f_8 &=& 2f^\pr _6\\
f_7^2 +tf_8f_7 +f_8 &=& 2f^\pr_7,
\end{array}
\right.
\eeq
where $ t=g_x(u,u)$,  $ f_i = f_i^F,\ i=6,7,8$
 and $ f^\prime_i$  denotes the first  derivative \\ of $ f_i$.\\
Then the first equation of the system  $(\ref{S})$ gives 
\beq\label{k1}
f_7(1+tf_6)=f_6
\eeq
and so $1+tf_6\neq 0,\quad \forall t\geq 0$, (otherwise
$1+tf_6= 0\quad \mbox{ would imply }\\ 
f_6 = 0\quad \mbox{ and } \quad tf_6=-1$, which is  absurd).
Hence (\ref{k1}) gives
\beq\label{kk1}
f_7 = \frac{f_6}{1+tf_6}.
\eeq
Furthermore the second equation of $(\ref{S})$ gives
\beq\label{k2}
f_8 = \frac{2f_6'-f_6^2}{1+tf_6}\quad .
\eeq
Next by using (\ref{kk1}), we obtain 
\beq\label{k3}
f_7^\pr =\frac{f_6^\pr-f_6^2}{(1+tf_6)^2}\quad ,
\eeq
and
\beq\label{k4}
1+tf_7 = \frac{1+2tf_6}{1+tf_6}\quad .
\eeq
By replacing (\ref{kk1}), (\ref{k2}), (\ref{k4})  and (\ref{k3})
in the $3^{rd}$ equation of the sytem $(\ref{S})$, we obtain
\beq\label{k5}
4tf_6f_6^\pr =-2f_6^2 +2tf_6^3,
\eeq
which implies
\beq\label{k6}
 f_6(t)&=&0\ ,\qquad \mbox{ or}\,   \\  
f_6^\pr(t) &=& -\frac{f_6(t)}{2t} +\frac{f_6^2(t)}{2}\, ;\\ 
\mbox{ for } t>0.    \nb
\eeq
So $f_6$ is a solution on the open set 
$I=\{t\in \, ]0\, ,\, +\infty[\; /\; f_6(t)\neq 0\}$  of 
the Bernouilli equation
\beq
y^\pr(t) &=& -\frac{y(t)}{2t} +\frac{y^2(t)}{2}\, .  
\label{k7}
\eeq
Besides, we have $f_6(0) = 0$. Indeed, if $0\in Adh(I)$ 
the adherence of $I$ in $\mbb{R}^+$, then by equation (\ref{k7}), we have 
\beq 
f_6(0)& =&\lim\limits_{\stackrel{t \to 0}{t\in I}} f_6(t)\nb \\ \nb 
 &=&\lim\limits_{\stackrel{t \to 0}{t\in I}}t[-2  f_6^\pr(t)+    f_6^2(t)]=0 .
 \eeq
But if $0\not\in Adh(I)$ then evidently, we have $f_6(0) =0$.

Thus the frontier $Fr(I)$ of $I$ is necessarily non empty, since $\mbb{R}^+$ is connected
and $f_6$ is smooth. In summary $f_6$ is a solution of the
equation
\beq
\left\{
\begin{array}{l}
y^\prime(t) = -\frac{y(t)}{2t} +y^2(t),\; \forall t\in I,\\
y_{|_{Fr(I)}} \equiv 0
\end{array}
\right.
\eeq
that has the unique solution $y\equiv 0$, so $f_6\equiv 0$.\\
Next by using (\ref{kk1}) and (\ref{k2}), we obtain  
$f_7=f_8 =0$, as stated.

                        $\hfill{\square}$
\bthm\label{TC}
Let $(M,g)$ be a Riemannian manifold and  $\ (TM,G)$  its \\
tangent bundle equipped with a $g$-natural metric G. Then $(TM,G)$ is flat\\
Riemannian if and only if
\begin{enumerate}
\item[i)] $(M,g)$ is flat, 
\item[ii)]
$\al_1(t)>0,\quad  \phi_1(t)>0,\quad
\al(t)> 0,\quad \phi(t)>0,\quad $
for all $t\in \mbb{R}^+$; 
\item[iii)] $\al_1+\al_3= constant>0,\quad \be_1+\be_3=0,\quad
2\al_2^\pr=\be_2,\quad $  
\item[iv)] $\al_1^\pr=\frac{\al_2\be_2}{\al_1+\al_3}\quad$
and $\quad \be_1 =\frac{\be_2(2\al_2+t\be2)}{\al_1+\al_3}$,
 \end{enumerate}
 where $\ \al_1^\pr\ $ and $\ \al_2^\pr $ are respectively the first derivatives
 of the functions \\
 $\ \al_1\ $ and $\ \al_2\ $.
  \ethm

\underline{{ \it Proof}}
\vspace{0.25cm}

Let us assume that $(TM,G)$ is flat Riemannian. 
By   Proposition \ref{P2} and   Proposition \ref{rp1},
we obtain the parts  $i)$ and $ii)$ of  Theorem \ref{TC}.

Next we  obtain $iii)$ from   Lemma \ref{L5}.

It remains to prove $iv)$. But according to Lemma \ref{L5} we have
\beq 
2\al_2^\pr &=&\be_2\quad \mbox{ and }\label{y1} \\ 
f_6&=&-\frac{\al_2}{\al}(\al_2^\pr +\frac{1}{2}\be_2)
+\al_1'\frac{(\al_1+\al_3)}{\al}=0. \label{y2}
\eeq
Then  by combining these identities, we obtain
\beq\label{y3}
\al_1^\pr = \frac{\al_2\be_2}{\al_1+\al_3}\ .
\eeq
Lemma \ref{L5} gives again 
\beq
f_7&=&(\be_1-\al_1^\pr)(\phi_1+\phi_3)-\be_2\phi_2=0\ ,\quad \mbox{ and }\\ \nb 
\be_1+\be_3&=&0, \nb
\eeq
then
\beq
\be_1&=& \al_1^\pr +\frac{\be_2\phi_2}{\al_1+\al_3}\\ \nb
 &=&\frac{\al_2\be_2}{\al_1+\al_3} +\frac{\be_2(\al_2+t\be_2)}{\al_1+\al_3}
 \quad \mbox{by}\quad (\ref{y3}) \\
 \nb
\be_1&=& \frac{\be_2(2\al_2+t\be_2)}{\al_1+\al_3}.
\eeq
So we prove $iv)$.

Conversely: \\
The part $ii)$ shows that $G$ is Riemannian. Next by combining the parts\\
$i)$ and  $iii)$ we obtain 
\beq\label{z1}
A=B=C=D=0\ .
\eeq
Furthermore by combining the parts $iii)$  and  $iv)$   
we obtain 
\beq
f_6^F&=&f_7^F=f_8^F=0,\label{z2}\\
f_6^E&=&f_7^E= \frac{\be_2}{\al_1+\al_3}\ ,
\;   f_8^E = \frac{2\be_2'}{\al_1+\al_3}\, . \label{z3}
\eeq
So (\ref{z2}) implies that $F=0$, and by considering (\ref{z1}) we obtain:\\
 $\forall\ (x,u)\in TM \mbox{ and } \forall\,  X,Y,Z\in T_xM$, 
 \beq
 \bar{R}\left(X^h,Y^h\right)Z^h &=& \bar{R}\left(X^h,Y^h\right)Z^v=0\\ \nb
 \bar{R}\left(X^h,Y^v\right)Z^h&=&\bar{R}\left(X^h,Y^v\right)Z^v = 
 \bar{R}\left(X^v,Y^v\right)Z^h = 0\; ,\nb
 \eeq
 where the lifts are taken at $(x,u)$.
Next  (\ref{z3}) implies
\beq
\bar{R}\left(X^v,Y^v\right)Z^v &=& h\{ d\left(E_{(Y,Z)}\right)_u(X)
                    -d\left(E_{(X,Z)}\right)_u(Y)\} \\ \nb
		 &=&\{ (f_7^E-f_6^E)g(Y,Z) 
		 + (f_8^E-2{f_6^E}')g(Y,u)g(Z,u)\}X\\ \nb
	  &&-\{ (f_7^E-f_6^E)g(X,Z)
	                   + (f_8^E-2{f_6^E}')g(X,u)g(Z,u)\}Z\\ \nb
		  &&+(2{f_7^E}'-f_8^E)\{g(Y,Z)g(X,u)-g(X,Z)g(Y,u)\}u\\ \nb
	                &=& 0\ .
\eeq
Finally $\bar{R}\equiv 0\ .\ $ 

   $\hfill{\square}$

\bibliographystyle{plain}

\end{document}